\documentclass[11pt]{article}
\usepackage[matrix,arrow,curve,cmtip]{xy}
\usepackage{amssymb}
\usepackage{latexsym}
\usepackage{theorem}
\usepackage[a4paper,width=140mm,height=220mm]{geometry}

\def\to{\mbox{$\xymatrix@1@C=5mm{\ar@{->}[r]&}$}}
\def\tto{\mbox{$\xymatrix@1@C=5mm{\ar@{=>}[r]&}$}}
\def\halfcirc{\begin{picture}(0,0)\put(0,3){\oval(2,2)[l]}\end{picture}}
\def\incl{\mbox{$\xymatrix@1@C=5mm{\ar@{->}[r]|<{\halfcirc}&}$}}
\def\bkar{\mbox{$\xymatrix@1@C=5mm{\ar@{->}[l]&}$}}
\def\distsign{\begin{picture}(0,0)\put(0,0){\circle{4}}\end{picture}}
\def\dist{\mbox{$\xymatrix@1@C=5mm{\ar@{->}[r]|{\distsign}&}$}}
\def\bkdist{\mbox{$\xymatrix@1@C=5mm{\ar@{->}[l]|{\distsign}&}$}}
\def\biar{\mbox{$\xymatrix@1@C=5mm{\ar@<1.5mm>[r]\ar@<-0.5mm>[r]&}$}}
\def\bidist{\mbox{$\xymatrix@1@C=5mm{\ar@<1.5mm>[r]|{\distsign}\ar@<-0.5mm>[r]|{\distsign}&}$}}
\def\adjar{\mbox{$\xymatrix@1@C=5mm{\ar@<1.5mm>@{<-}[r]\ar@<-0.5mm>[r]&}$}}
\def\adjdist{\mbox{$\xymatrix@1@C=5mm{\ar@<1.5mm>@{<-}[r]|{\distsign}\ar@<-0.5mm>[r]|{\distsign}&}$}}
\def\iso{\mbox{$\xymatrix@1@C=6mm{\ar@{->}[r]^{\sim}&}$}}
\def\doubiso{\mbox{$\xymatrix@1@C=6mm{\ar@{<->}[r]^{\sim}&}$}}
\def\doubar{\mbox{$\xymatrix@1@C=6mm{\ar@{<->}[r]&}$}}

\def\endoar#1#2{\mbox{\xymatrix{{#1}\ar@(u,r)|{#2}}}}

\newtheorem{theorem}{Theorem}[section]
\newtheorem{lemma}[theorem]{Lemma}
\newtheorem{definition}[theorem]{Definition} 
\newtheorem{proposition}[theorem]{Proposition}
\newtheorem{corollary}[theorem]{Corollary}
{\theorembodyfont{\upshape}\newtheorem{example}[theorem]{Example}}
\newcommand{\proof}{\noindent {\it Proof\ }: }
\def\endofproof{$\mbox{ }\hfill\Box$\par\vspace{1.8mm}\noindent}

\def\wh#1{\widehat{#1}}
\def\wt#1{\widetilde{#1}}
\def\TRSDist{{\sf TRSDist}}
\def\TRSCat{{\sf TRSCat}}

\def\Idm{{\sf Idm}}

\def\R{{\cal R}}
\def\RSDist{{\sf RSDist}}
\def\RSCat{{\sf RSCat}}

\def\wh#1{\widehat{#1}}
\def\ol#1{\overline{#1}}
\def\+{^{\dagger}}

\def\norm{_{\sf n}}

\def\Ord{{\sf Ord}}
\def\Idl{{\sf Idl}}

\def\Sh{{\sf Sh}}
\def\Rel{{\sf Rel}}

\def\:{\colon}
\def\1{{\bf 1}}

\def\2{{\bf 2}}
\def\3{{\bf 3}}
\def\inv{^{-1}}

\def\Set{{\sf Set}}
\def\QUANT{{\sf QUANT}}

\def\cc{_{\sf cc}}

\def\dom{{\sf dom}}

\def\Set{{\sf Set}}

\def\Sup{{\sf Sup}}
\def\Cat{{\sf Cat}}
\def\Matr{{\sf Matr}}

\def\Dist{{\sf Dist}}

\def\Cat{{\sf Cat}}

\def\SCat{{\sf SCat}}
\def\Map{{\sf Map}}

\def\Q{{\cal Q}}

\def\P{{\cal P}}
\def\V{{\cal V}}

\def\lim{\mathop{\rm lim}}
\def\bbA{\mathbb{A}}
\def\bbB{\mathbb{B}}
\def\bbC{\mathbb{C}}
\def\bbD{\mathbb{D}}

\def\tensor{\otimes}
\def\<{\langle}
\def\>{\rangle}

\title{Categorical structures enriched in a quantaloid: \\ orders and ideals over a base quantaloid}
\author{Isar Stubbe\footnote{D\'epartement de Math\'ematique, Universit\'e de Louvain, Chemin du Cyclotron 2,
1348 Louvain-la-Neuve (Belgique), {\tt
i.stubbe@math.ucl.ac.be}.}}
\date{August 7, 2004}

\begin{document}

\maketitle

\begin{quote}{\small

{\bf Abstract.} Applying (enriched) categorical structures we define the notion of ordered sheaf on a quantaloid $\Q$, which we call `$\Q$-order'. This requires a theory of semicategories enriched in the quantaloid $\Q$, that admit a suitable Cauchy completion. There is a quantaloid $\Idl(\Q)$ of $\Q$-orders and ideal relations, and a locally ordered category $\Ord(\Q)$ of $\Q$-orders and monotone maps; actually, $\Ord(\Q)=\Map(\Idl(\Q))$. In particular is $\Ord(\Omega)$, with $\Omega$ a locale, the category of ordered objects in the topos of sheaves on $\Omega$. In general $\Q$-orders can equivalently be described as Cauchy complete categories enriched in the split-idempotent completion of $\Q$. Applied to a locale $\Omega$ this generalizes and unifies previous treatments of (ordered) sheaves on $\Omega$ in terms of $\Omega$-enriched structures.

}\end{quote}

\section{Introduction}\label{A}

An ``ordered set'' is usually thought of as a set equipped with a reflexive and transitive relation; that is to say, it is an ordered object in $\Set$. But one can also treat an order $(A,\leq)$ by means of the classifying map for its order relation, say $[\cdot\leq\cdot]\:A\times A\to\2$, where now $\2$ is the object of truth values. This takes us into the realm of enriched categorical structures, for the reflexivity and transitivity axioms on the order relation translate into unit-inequalities and composition-inequalities for the enrichment $[\cdot\leq\cdot]$ of $A$ over $\2$. So order theory is then a matter of applied (enriched) categorical structures.
\par
More generally, an ``ordered sheaf on a locale $\Omega$'' is an ordered object in the topos $\Sh(\Omega)$ of sheaves on the locale. Here too one may attempt at describing such an $\Omega$-order $(A,\leq)$ in terms of enriched categorical structures. There are two approaches: [Walters, 1981] (implicitly) treats such $\Omega$-orders as {\em categories} enriched in $\Rel(\Omega)$; whereas [Borceux and Cruciani, 1998] prefer to work with {\em semicategories} enriched in $\Omega$. The first option has the advantage that it speaks of categories enriched in a quantaloid, a clear and transparent theory that may be developed along the lines of the well-known theory of $\V$-enriched categories; but it has the disadvantage that its base quantaloid is $\Rel(\Omega)$ and not $\Omega$ itself, so that the r\^ole of the truth values is slightly obscured. The second option prominently keeps the locale $\Omega$ as base for enrichment; but the price to pay is that one has to work with enriched semicategories, i.e.~``categories without units''. These two approaches must be equivalent of course, but an explicit comparison was not yet provided in the literature.
\par
The subject of this paper is to go still a bit further, and to consider ``ordered sheaves on a quantaloid $\Q$''. A quantaloid is a $\Sup$-enriched category, and a quantale is a one-object quantaloid i.e.~a monoid in $\Sup$; so a locale is a very particular quantaloid. Unfortunately there is no ``topos of sheaves on a quantaloid $\Q$'', so we cannot define a ``$\Q$-order'' as ordered object in such a topos. However, we can still consider categorical structures enriched in a quantaloid, and using those we will define what we believe to be the correct notion of ``$\Q$-order''.
\par
Our option has been, in a first instance, to (clarify and) generalize the work of [Borceux and Cruciani, 1998] from locale-enrichment to quantaloid-enrichment. So we work with $\Q$-enriched semicategories, building further on the notions introduced in [Stubbe, 2004b]. We analyze what it takes for a $\Q$-semicategory to admit a well-behaved Cauchy completion (for indeed this is non-trivial) and come up with a notion of ``totally regular $\Q$-semicategory''. In our opinion, the Cauchy complete totally regular $\Q$-semicategories are then precisely the sought-after ``$\Q$-orders''. We define a (locally ordered) category $\Ord(\Q)$ of ``$\Q$-orders and monotone maps'' and a quantaloid $\Idl(\Q)$ of ``$\Q$-orders and ideal relations'', and we have that $\Ord(\Q)=\Map(\Idl(\Q))$. Taking a base locale $\Omega$ gives back the description of $\Omega$-orders as put forward by [Borceux and Cruciani, 1998].
\par
But we also show how in general $\Q$-orders can be seen as Cauchy complete categories enriched in the split-idempotent completion of $\Q$. Then working with a locale $\Omega$ gives back the description of $\Omega$-orders (implicitly) given by [Walters, 1981]. As a result we then have a generalization of, and an explicit comparison between, Lovanists' and Sydneysiders' previous work on the subject of ordered objects in $\Sh(\Omega)$ in terms of enriched categorical structures.

\section{Background}\label{B}

Throughout this paper we will rely heavily on the theory of categorical structures enriched in a base quantaloid. To make this paper reasonably self-contained we recall some definitions and notations; for details we refer to [Stubbe, 2004a, 2004b] (where more references can be found).

\subsection*{Quantaloids}

A {\em quantaloid} is a category enriched in the symmetric monoidal closed category $\Sup$ of complete lattices and morphisms that preserve arbitrary suprema. That is to say, such is a category whose hom-objects are complete lattices, and composition of morphisms distributes on both sides over arbitrary suprema of morphisms. A {\em homomorphism} of quantaloids is a $\Sup$-functor between $\Sup$-enriched categories: it is a functor that preserves arbitrary suprema of morphisms. Any quantaloid is in particular a closed bicategory: we denote
$$-\circ f\dashv \{f,-\}\mbox{ and }f\circ-\dashv[f,-]$$
for the adjoints to composition with some morphisms $f\:A\to B$ in a quantaloid. Many calculations in a quantaloid depend on its closedness.
\par
A {\em quantale} is, by definition, a monoid in $\Sup$; in other words, it is a quantaloid with only one object. Since a {\em locale} is precisely a quantale with multiplication given by infima, it too may be viewed as a quantaloid (with only one object). Therefore the theory of quantaloid-enriched categorical structures applies to quantales and locales.
\par
From now on, let $\Q$ denote a {\em small} quantaloid.

\subsection*{Quantaloid-enriched structures}

A {\em $\Q$-typed set $X$} is an object of the slice category $\Set/\Q_0$ of sets over the object-set of $\Q$. In other terms, such is a set $X$ to every element of which is associated an object of $\Q$: for every $x\in X$ there is a $tx$ in $\Q$ (which is called the {\em type} of $x$ in $\Q$). The notation with a ``$t$'' for the types of elements in a $\Q$-typed set is generic: even for two different $\Q$-typed sets $X$ and $Y$, the type of $x\in X$ is written $tx$ and that of $y\in Y$ is $ty$.
\par
A {\em $\Q$-enriched category $\bbA$} is then a $\Q$-typed set $\bbA_0$ (of ``objects'') together with, for every $(a,a')\in\bbA_0\times\bbA_0$, a $\Q$-arrow $\bbA(a',a)\:ta\to ta'$ (``the hom-arrow from $a$ to $a'$'') satisfying
$$\bbA(a'',a')\circ\bbA(a',a)\leq\bbA(a'',a)\mbox{ and }1_{ta}\leq\bbA(a,a)$$
for all $a,a',a''\in\bbA_0$. A {\em distributor} $\Phi\:\bbA\dist\bbB$ between two $\Q$-categories is a matrix of $\Q$-arrows $\Phi(b,a)\:ta\to tb$, one for each $(a,b)\in\bbA_0\times\bbB_0$, satisfying
$$\bbB(b',b)\circ\Phi(b,a)\leq\Phi(b',a)\mbox{ and }\Phi(b,a)\circ\bbA(a,a')\leq\Phi(b,a')$$
for every $a,a'\in\bbA_0$ and $b,b'\in\bbB_0$. The composite of two distributors $\Phi\:\bbA\dist\bbB$ and $\Psi\:\bbB\dist\bbC$ is written $\Psi\tensor\Phi\:\bbA\dist\bbC$ and is defined by
\begin{equation}\label{1}
(\Psi\tensor\Phi)(c,a)=\bigvee_{b\in\bbB_0}\Psi(c,b)\circ\Phi(b,a).
\end{equation}
The identity distributor on a $\Q$-category $\bbA$ is $\bbA$ itself: $\bbA(-,-)\:\bbA\dist\bbA$. Parallel distributors are ordered ``elementwise''; in particular is the supremum of certain $\Phi_i\:\bbA\dist\bbB$ ($i\in I$) given by
\begin{equation}\label{2}
(\bigvee_i\Phi_i)(b,a)=\bigvee_i\Phi_i(b,a).
\end{equation}
It may thus be verified that $\Q$-categories and distributors form a (large) quantaloid $\Dist(\Q)$. A {\em functor} $F\:\bbA\to\bbB$ on the other hand is an object mapping $\bbA_0\to\bbB_0\:a\mapsto Fa$ satisfying
$$t(Fa)=ta\mbox{ and }\bbA(a',a)\leq\bbB(Fa',Fa)$$
for all $a,a'\in\bbA_0$. $\Q$-categories and functors form a category $\Cat(\Q)$ in the obvious way. Crucially, every functor $F\:\bbA\to\bbB$ induces an adjoint pair of distributors: denote $\bbB(-,F-)\:\bbA\dist\bbB$ for the matrix of the $\Q$-arrows $\bbB(b,Fa)\:ta\to tb$, and likewise $\bbB(F-,-)\:\bbB\dist\bbA$ for the matrix of the $\bbB(Fa,b)\:tb\to ta$, then $\bbB(-,F-)\dashv\bbB(F-,-)$ in $\Dist(\Q)$. There is a functor 
$$\Cat(\Q)\to\Dist(\Q)\:\Big(F\:\bbA\to\bbB\Big)\mapsto\Big(\bbB(-,F-)\:\bbA\dist\bbB\Big)$$
which in turn induces a local ordering on $\Cat(\Q)$: for parallel functors $F,G\:\bbA\biar\bbB$ we put $F\leq G$ whenever $\bbB(-,F-)\leq\bbB(-,G-)$. So the functor above becomes a faithful 2-functor from a locally ordered category into a quantaloid.
\par
Next we consider ``categories without units''. A {\em $\Q$-enriched semicategory} is a $\Q$-typed set $\bbA_0$ of ``objects'' together with, for every $(a,a')\in\bbA_0\times\bbA_0$, a ``hom-arrow'' $\bbA(a',a)\:ta\to ta'$ in $\Q$ satisfying only
$$\bbA(a'',a')\circ\bbA(a',a)\leq\bbA(a'',a)$$
for all $a,a',a''\in\bbA_0$. A {\em semifunctor} $F\:\bbA\to\bbB$ between $\Q$-semicategories is defined in the same manner as a functor between $\Q$-categories, and there is a category $\SCat(\Q)$ of $\Q$-semicategories and semifunctors. The definition of {\em semidistributor} $\Phi\:\bbA\dist\bbB$ between two $\Q$-semicategories too is an exact copy of that of distributor between $\Q$-categories, but $\Q$-semicategories and semidistributors {\em do not form a quantaloid}, and {\em a fortiori} semifunctors between $\Q$-semicategories {\em do not induce adjoint pairs of semidistributors}. Because of this, the theory of $\Q$-semicategories doesn't allow for much developments: one has to add ``regularity conditions'' on those $\Q$-semicategories, semidistributors and semifunctors to make an interesting theory (with interesting examples). So the important generalization of the notion of $\Q$-category is the following: a {\em regular $\Q$-semicategory $\bbA$} is a $\Q$-typed set $\bbA_0$ of ``objects'' together with a ``hom-arrow'' $\bbA(a',a)\:ta\to tb$ in $\Q$ for every $(a,a')\in\bbA_0\times\bbA_0$ satisfying
$$\bigvee_{a'\in\bbA_0}\bbA(a'',a')\circ\bbA(a',a)=\bbA(a'',a)$$
for all $a,a''\in\bbA_0$. Every $\Q$-category is a regular $\Q$-semicategory, but the converse is not true. A {\em regular semidistributor} $\Phi\:\bbA\dist\bbB$ between regular $\Q$-semicategories is a matrix of $\Q$-arrows $\Phi(b,a)\:ta\to tb$ (one for every $(a,b)\in\bbA_0\times\bbB_0$) satisfying
$$\bigvee_{b\in\bbB_0}\bbB(b',b)\circ\Phi(b,a)=\Phi(b',a)\mbox{ and }\bigvee_{a\in\bbA_0}\Phi(b,a)\circ\bbA(a,a')=\Phi(b,a')$$
for every $a,a'\in\bbA_0$ and $b,b'\in\bbB_0$. Reusing the formula in (\ref{1}) the composition $\Psi\tensor\Phi\:\bbA\dist\bbB$ of two regular semidistributors $\Phi\:\bbA\dist\bbB$ and $\Psi\:\bbB\dist\bbC$ between regular $\Q$-semicategories is defined, and reusing (\ref{2}) gives the supremum of parallel regular semidistributors. Also, the identity regular semidistributor on a regular $\Q$-semicategory $\bbA$ is $\bbA(-,-)\:\bbA\dist\bbA$ itself. So there is a (large) quantaloid $\RSDist(\Q)$ of regular $\Q$-semicategories and regular semidistributors, of which $\Dist(\Q)$ is a full subquantaloid in the obvious way. Further, given two regular $\Q$-semicategories $\bbA$ and $\bbB$, a {\em regular semifunctor} $F\:\bbA\to\bbB$ is an object mapping $\bbA_0\to\bbB_0\:a\mapsto Fa$ that satisfies
$$t(Fa)=ta\mbox{, }\bbA(a',a)\leq\bbB(Fa',Fa),$$
$$\bigvee_{a'\in\bbA_0}\bbB(b,Fa')\circ\bbA(a',a)=\bbB(b,Fa)\mbox{ and }\bigvee_{a'\in\bbA_0}\bbA(a,a')\circ\bbB(Fa',b)=\bbB(Fa,b)$$
for all $a,a'\in\bbA_0$ and $b\in\bbB_0$. The point of the latter two conditions is precisely to assure that both matrices $\bbB(-,F-)$ and $\bbB(F-,-)$ are regular semidistributors; in fact, mimicking the case of $\Q$-categories and functors, here too $\bbB(-,F-)\dashv\bbB(F-,-)$ in the quantaloid $\RSDist(\Q)$. Regular $\Q$-semicategories and regular semifunctors form a category $\RSCat(\Q)$ in the obvious way, and $\Cat(\Q)$ is a full subcategory; the action
$$\RSCat(\Q)\to\RSDist(\Q)\:\Big(F\:\bbA\to\bbB\Big)\mapsto\Big(\bbB(-,F-)\:\bbA\dist\bbB\Big)$$
is functorial, and gives rise to a local ordering on $\RSCat(\Q)$: for regular semifunctors $F,G\:\bbA\biar\bbB$ we put $F\leq G$ precisely when $\bbB(-,F-)\leq\bbB(-,G-)$. This makes the action above a faithful 2-functor from a locally ordered category into a quantaloid.
\par
The following commutative diagram of 2-categories and 2-functors summarizes:
$$\xymatrix@=15mm{
\Cat(\Q)\ar[r]\ar[d] & \RSCat(\Q)\ar[d] \\
\Dist(\Q)\ar[r] & \RSDist(\Q)}$$
The horizontal arrows are full embeddings saying that the theory of regular $\Q$-semicategories is more general than that of $\Q$-categories; and the vertical arrows are faithful 2-functors saying that every (regular semi)functor between (regular semi)categories induces a left adjoint (regular semi)distributor.

\subsection*{Cauchy complete $\Q$-categories}

A (necessarily left adjoint) distributor $\Phi\:\bbA\dist\bbB$ between $\Q$-categories is said to {\em converge} when there exists a (necessarily essentially unique) functor $F\:\bbA\to\bbB$ such that $\bbB(-,F-)=\Phi$. Using the term {\em Cauchy distributor} synonymously for left adjoint distributor, if all Cauchy distributors into a $\Q$-category $\bbB$ converge, then $\bbB$ is by definition {\em Cauchy complete}. Thus it follows trivially that $\Cat\cc(\Q)\simeq\Map(\Dist\cc(\Q))$
when we consider the obvious full subcategories determined by the Cauchy complete $\Q$-categories. 
\par
For every $\Q$-category $\bbB$ there exists a Cauchy complete $\Q$-category $\bbB\cc$ such that $\bbB$ and $\bbB\cc$ are isomorphic objects of $\Dist(\Q)$. So $\Dist\cc(\Q)=\Dist(\Q)$, and in particular $\Cat\cc(\Q)\simeq\Map(\Dist(\Q))$.
The construction of $\bbB\cc$ goes as follows: its objects are the left adjoint distributors
into $\bbB$ whose domain is a one-object $\Q$-category with as (single) hom-arrow an identity in $\Q$, like\footnote{For an object $X$ of $\Q$ we denote $*_X$ for the $\Q$-category with one object and hom-arrow $1_X$.} $\phi\:*_X\dist\bbB$, and the type of such an object is $t\phi=X$; for two such objects, say $\phi\:*_X\dist\bbB$ and $\psi\:*_Y\dist\bbB$, the hom-arrow $\bbB\cc(\psi,\phi)\:X\to Y$ is (the single element of) the lifting $[\psi,\phi]$ in $\Dist(\Q)$.
\par
Every object $b$ of a $\Q$-category $\bbB$ can be ``pointed at'' by a canonical constant functor; to wit, 
$$\Delta b\:*_{tb}\to\bbB\:*\mapsto b.$$
As for any functor, there is an associated adjoint pair of distributors:
$$\xymatrix@=15mm{
\ast_{tb}\ar@{}[r]|{\perp}\ar@/^3mm/@<1mm>[r]|{\distsign}^{\bbB(-,b)} & \bbB\ar@/^3mm/@<1mm>[l]|{\distsign}^{\bbB(b,-)}}.$$
Therefore there is an embedding 
$$k_{\bbB}\:\bbB\to\bbB\cc\:b\mapsto\Big(\bbB(-,b)\:*_{tb}\dist\bbB\Big)$$
which, by the Yoneda lemma for $\Q$-categories, is fully faithful. (It is actually the distributor induced by the functor $k_{\bbB}$ that gives the aforementioned isomorphism $\bbB\cong\bbB\cc$ in $\Dist(\Q)$.) Therefore $\bbB\cc$ merits to be called the {\em Cauchy completion of $\bbB$}; and $\bbB$ may be thought of as full subcategory of $\bbB\cc$ determined by the ``constant'' objects.

\section{Total regularity}\label{C}

Awkward as it may seem, the objects of a regular $\Q$-semicategory $\bbA$ need not be ``constant'': they can in general not be ``pointed at'' by a constant (regular semi)functor. This obstructs a well-behaved Cauchy completion of $\bbA$. So we need the following definition.
\begin{definition}\label{3}
An object $a$ of a regular $\Q$-semicategory $\bbA$ is {\em stable} when there exists a regular semifunctor from a one-object regular $\Q$-semicategory into $\bbA$ whose value is $a$. And the regular $\Q$-semicategory $\bbA$ is {\em totally regular} when all of its objects are stable.
\end{definition}
Not every regular $\Q$-enriched semicategory is totally regular, and not every totally regular $\Q$-semicategory is a $\Q$-category. On the other hand, any $\Q$-category is totally regular. Denoting $\TRSCat(\Q)$ and $\TRSDist(\Q)$ for the full substructures of $\RSCat(\Q)$ and $\RSDist(\Q)$ determined by the totally regular $\Q$-semicategories, we have a diagram of 2-categories and 2-functors
$$\xymatrix@=15mm{
\Cat\ar[r]\ar[d] & \TRSCat(\Q)\ar[r]\ar[d] & \RSCat(\Q)\ar[d] \\
\Dist(\Q)\ar[r] & \TRSDist(\Q)\ar[r] & \RSDist(\Q)}$$
in which all horizontal arrows are full embeddings, and the vertical ones are faithful 2-functors saying that every (regular semi)functor between ((totally) regular semi)categories induces an adjoint pair of (regular semi)distributors. In the next sections we will argue that {\em the totally regular $\Q$-semicategories are precisely those semicategories that admit a suitable Cauchy completion}. But before doing so we will analyze the definition of `total regularity'.
\par
The one-object regular $\Q$-semicategories correspond precisely to the
idempotents\footnote{More generally, regular $\Q$-semicategories are idempotents in $\Matr(\Q)$, the quantaloid of matrices with elements in $\Q$.} in $\Q$: any idempotent $\Q$-arrow $e$ determines a regular $\Q$-semicategory $*_e$ with one object which is of type $\dom(e)$, and whose single hom-arrow is precisely $e$; conversely, the single hom-arrow of a one-object regular $\Q$-semicategory is an idempotent $\Q$-arrow.
\begin{lemma}\label{4}
For an object $a$ of a regular $\Q$-semicategory $\bbA$, the following
are equivalent:
\begin{enumerate}
\item\label{5} $a$ is stable;
\item\label{6} there exists an idempotent $\Q$-arrow $e\:ta\to ta$ such that $e\leq\bbA(a,a)$ and, for
all
$a'\in\bbA_0$, $\bbA(a',a)\circ e=\bbA(a',a)$ and $e\circ\bbA(a,a')=\bbA(a,a')$;
\item\label{7} for all $a'\in\bbA_0$, $\bbA(a',a)\circ \bbA(a,a)=\bbA(a',a)$ and $\bbA(a,a)\circ\bbA(a,a')=\bbA(a,a')$.
\end{enumerate}
\end{lemma}
\proof The equivalence of \ref{5} and \ref{6} is straightforward from the definitions. 
Assuming \ref{6}, in the regular $\Q$-semicategory $\bbA$ we have
\begin{eqnarray*}
\bbA(a',a)\circ\bbA(a,a)
 & \leq & \bigvee_{x\in\bbA_0}\bbA(a',x)\circ\bbA(x,a) \\
 & = & \bbA(a',a) \\
 & = & \bbA(a',a)\circ e \\
 & \leq & \bbA(a',a)\circ\bbA(a,a)
\end{eqnarray*}
which proves that $\bbA(a',a)\circ\bbA(a,a)=\bbA(a',a)$. Similar for the other equation in \ref{7}. Conversely, putting $a=a'$ in \ref{7} shows that $\bbA(a,a)\:ta\to ta$ is an idempotent in $\Q$ satisfying the conditions in \ref{6}.
\endofproof
This lemma says that an object $a$ of a regular $\Q$-semicategory is stable if and only if the $\Q$-arrow
$\bbA(a,a)$ is idempotent and plays the r\^ole of ``identity at $a$''.
So there is a ``canonical'' regular semifunctor pointing at such a stable object $a$: 
$$\Gamma a\:*_{\bbA(a,a)}\to\bbA\:*\mapsto a.$$ 
Like any regular semifunctor, $\Gamma a$ induces an adjoint pair of regular semidistributors:
\begin{equation}\label{7.1}
\ast_{\bbA(a,a)}\xymatrix@=15mm{
\ar@/^3mm/@<1mm>[r]^{\bbA(-,a)}|{\distsign}\ar@{}[r]|{\perp} & \ar@/^3mm/@<1mm>[l]^{\bbA(a,-)}|{\distsign}}\bbA.
\end{equation}
This will be important later on: it will allow us to identify the objects of a totally regular $\Q$-semicategory $\bbA$ with the ``constant'' objects of its Cauchy completion $\bbA\cc$ (see \ref{14} and further).
\par
As a consequence of \ref{4} we can now give an elementary characterization of `totally regular $\Q$-semicategory'.
\begin{proposition}\label{8}
A totally regular $\Q$-semicategory $\bbA$ is a $\Q$-typed set $\bbA_0$ (of ``objects'') and, for every pair $(a,a')\in\bbA_0\times\bbA_0$, a $\Q$-arrow $\bbA(a',a)\:ta\to ta'$ 
(the ``hom-arrow from $a$ to $a'$'') satisfying, for all $a,a',a''\in\bbA_0$,
$$\bbA(a'',a')\circ\bbA(a',a)\leq\bbA(a'',a),$$
$$\bbA(a,a)\circ\bbA(a,a')=\bbA(a,a')\mbox{ and }
\bbA(a',a)\circ\bbA(a,a)=\bbA(a',a).$$
\end{proposition}
\proof
An $\bbA$ as in the statement of this proposition is a $\Q$-semicategory that in particular satisfies, for any $a,a'\in\bbA_0$,
$$\bbA(a',a)
\geq\bigvee_{x\in\bbA_0}\bbA(a',x)\circ\bbA(x,a)
\geq\bbA(a',a')\circ\bbA(a',a)=\bbA(a',a).$$
So $\bbA$ is actually a regular $\Q$-semicategory. But then \ref{4} implies that $\bbA$ is totally regular. Conversely it is -- due to \ref{4} -- trivial that a totally regular $\Q$-semicategory $\bbA$ verifies the (in)equations in the statement of this proposition.
\endofproof
Intuitively, a totally regular $\Q$-semicategory is a $\Q$-enriched semicategory (with composition law $\bbA(a'',a')\circ\bbA(a',a)\leq\bbA(a'',a)$), without global identities (no $1_{ta}\leq\bbA(a,a)$ in general) but with local identities ($\bbA(a,a)$ being an idempotent playing the r\^ole of ``identity at $a$''). 
\par
It is useful to also give elementary characterizations of `regular semidistributor' and `regular semifunctor' between totally regular $\Q$-semicategories.
\begin{proposition}\label{9}
Let $\bbA$ and $\bbB$ be two totally regular $\Q$-semicategories. A regular semidistributor $\Phi\:\bbA\dist\bbB$ is a matrix of $\Q$-arrows $\Phi(b,a)\:ta\to tb$, one for each couple $(a,b)\in\bbA_0\times\bbB_0$, satisfying for all $a,a'\in\bbA_0$ and $b,b'\in\bbB_0$,
$$\bbB(b',b)\circ\Phi(b,a)\leq\Phi(b',a)\mbox{, }
\Phi(b,a)\circ\bbA(a,a')\leq\Phi(b,a'),$$
$$\bbB(b,b)\circ\Phi(b,a)=\Phi(b,a)
\mbox{ and }\Phi(b,a)\circ\bbA(a,a)=\Phi(b,a).$$
\end{proposition}
\proof
We must show that a semidistributor $\Phi\:\bbA\dist\bbB$ between totally regular $\Q$-semicategories is regular if and only if it satisfies $$\bbB(b,b)\circ\Phi(b,a)=\Phi(b,a)
\mbox{ and }\Phi(b,a)\circ\bbA(a,a)=\Phi(b,a).$$
The calculations are analogous to those in the proof of \ref{8}.
\endofproof
\begin{proposition}\label{10}
Let $\bbA$ and $\bbB$ be two totally regular $\Q$-semicategories. A regular semifunctor $F\:\bbA\dist\bbB$ is an object mapping $\bbA_0\to\bbB_0\:a\mapsto Fa$ such that, for all $a,a'\in\bbA_0$ and $b\in\bbB_0$,
$$t(Fa)=ta\mbox{, }\bbA(a',a)\leq\bbB(Fa',Fa),$$
$$\bbB(b,Fa)\circ\bbA(a,a)=\bbB(b,Fa)\mbox{ and }
\bbA(a,a)\circ\bbB(Fa,b)=\bbB(Fa,b).$$
\end{proposition}
\proof
We must show that, for a semifunctor $F\:\bbA\to\bbB$ between (totally) regular $\Q$-semi\-cat\-e\-go\-ries, the semidistributors $\bbB(-,F-)\:\bbA\dist\bbB$ and $\bbB(F-,-)\:\bbB\dist\bbA$ are regular if and only if
$$\bbB(b,Fa)\circ\bbA(a,a)=\bbB(b,Fa)\mbox{ and }
\bbA(a,a)\circ\bbB(Fa,b)=\bbB(Fa,b).$$
This follows by application of \ref{9}.
\endofproof
\par
To end this section we make a helpful observation.
\begin{corollary}\label{11}
Consider two $\Q$-semicategories $\bbA$ and $\bbB$, and a ``fully faithful semifunctor'' $F\:\bbA\to\bbB$, i.e.~an object mapping $\bbA_0\to\bbB_0\:a\mapsto Fa$ that satisfies, for all $a,a'\in\bbA_0$, 
$$t(Fa)=ta\mbox{ and }\bbA(a',a)=\bbB(Fa',Fa).$$ 
If $\bbB$ is totally regular then so is $\bbA$, and then $F$ is a regular semifunctor.
\end{corollary}
\proof
We verify the conditions in \ref{8} for $\bbA$: for $a,a'\in\bbA_0$,
$$\bbA(a',a)\circ\bbA(a,a)=\bbB(Fa',Fa)\circ\bbB(Fa,Fa)=\bbB(Fa',Fa)=\bbA(a',a)$$
and likewise for $\bbA(a',a')\circ\bbA(a',a)=\bbA(a',a)$, so $\bbA$ is totally regular. With a similar trick one verifies that $F$ meets the conditions in \ref{9}.
\endofproof
In particular does it follow that every full subgraph of a totally regular $\Q$-semi\-cat\-e\-go\-ry is again a totally regular $\Q$-semicategory. For a $\Q$-category it is evident that any full subgraph is a $\Q$-category too. But a full subgraph of a merely regular $\Q$-semicategory is not necessarily regular!

\section{Cauchy completion}\label{D}

\subsection*{Cauchy complete totally regular semicategories}

For totally regular $\Q$-semicategories we will mimic a lot of the terminology and technology that is known for $\Q$-categories. So a left adjoint regular semidistributor is also called a {\em Cauchy regular semidistributor}; usually  we will denote the right adjoint of some $\Phi\:\bbA\dist\bbB$ in $\TRSDist(\Q)$ as $\Phi^*\:\bbB\dist\bbA$. Such a Cauchy regular semidistributor $\Phi\:\bbA\dist\bbB$ is said to {\em converge} when there exists a regular semifunctor\footnote{Note: a left adjoint $\Phi\:\bbA\dist\bbB$ in $\TRSDist(\Q)$ converges if and only if there exists a semifunctor $F\:\bbA\to\bbB$ such that {\em both} $\Phi=\bbB(-,F-)$ and $\Phi^*=\bbB(F-,-)$---in which case $F$ {\em proves} to be regular.} $F\:\bbA\to\bbB$ such that $\Phi=\bbB(-,F-)$ (equivalently, $\Phi^*=\bbB(F-,-)$). When $\Phi\:\bbA\dist\bbB$ converges to both $F\:\bbA\to\bbB$ and $G\:\bbA\to\bbB$, then $F$ and $G$ are equivalent in the order $\TRSCat(\Q)(\bbA,\bbB)$; so the notion of convergence is well-defined.
\begin{definition}\label{12}
A {\em Cauchy complete totally regular $\Q$-semicategory} is one for which all Cauchy regular semidistributors into it, converge. 
\end{definition}
In other words, an object $\bbB$ of $\TRSCat(\Q)$ is Cauchy complete if and only if, for every other object $\bbA$, 
$$\TRSCat(\Q)(\bbA,\bbB)\to\Map(\TRSDist(\Q))(\bbA,\bbB)\:F\mapsto\bbB(-,F-)$$
is surjective---in which case it is an equivalence of orders, since it is always ``essentially injective''.
\par
The next lemma asserts that, when verifying the Cauchy completeness of some totally regular $\Q$-semicategory $\bbB$, it suffices to verify the convergence of the left adjoint regular semidistributors into $\bbB$ {\em whose domain has only one object}. Such a left adjoint $\phi\:*_{e}\dist\bbB$ in $\TRSDist(\Q)$ (with $e$ an idempotent in $\Q$) converges if and only if there exists an object $b\in\bbB_0$ with type $tb=\dom(e)$ such that $e\leq\bbB(b,b)$, $\bbB(y,b)=\phi(y)$ and $\bbB(b,y)=\phi^*(y)$ for all $y\in\bbB_0$. (That is to say, these conditions are equivalent to saying that $*_e\to\bbB\:*\mapsto b$ is a regular semifunctor to which $\phi$ converges---see the previous footnote.)
\begin{lemma}\label{13}
A totally regular $\Q$-semicategory $\bbB$ is Cauchy complete if and only if, for every idempotent arrow $e$ in $\Q$, every Cauchy regular semidistributor $\phi\:*_{e}\dist\bbB$ converges.
\end{lemma}
\proof One implication is trivial. For the other, let $\bbA$ be a totally regular $\Q$-semicategory and $\Phi\:\bbA\dist\bbB$ a regular semidistributor with right adjoint $\Phi^*$. For every $a\in\bbA_0$ we may consider the composition of adjunctions in $\TRSDist(\Q)$ like so:
$$\ast_{\bbA(a,a)}\xymatrix@=15mm{
\ar@/^3mm/@<1mm>[r]^{\bbA(-,a)}|{\distsign}\ar@{}[r]|{\perp} & \bbA\ar@/^3mm/@<1mm>[l]^{\bbA(a,-)}|{\distsign}\ar@/^3mm/@<1mm>[r]^{\Phi}|{\distsign}\ar@{}[r]|{\perp} & \ar@/^3mm/@<1mm>[l]^{\Phi^*}|{\distsign}}\bbB.$$
By regularity, $\Phi\tensor\bbA(-,a)=\Phi(-,a)$ and $\bbA(a,-)\tensor\Phi^*=\Phi^*(a,-)$, and by assumption every such $\Phi(-,a)\:*_{\bbA(a,a)}\dist\bbB$ converges: there exist objects
$\{Fa\in\bbB_0\mid a\in\bbA_0\}$
with types $t(Fa)=\dom(\bbA(a,a))=ta$ such that $\bbA(a,a)\leq\bbB(Fa,Fa)$, $\bbB(y,Fa)=\Phi(y,a)$ and $\bbB(Fa,y)=\Phi^*(a,y)$.
Consider now the type preserving object mapping $\bbA_0\to\bbB_0\:a\mapsto Fa$: we verify that it satisfies all the conditions in \ref{10}. For $a,a'\in\bbA_0$ use the adjunction $\Phi\dashv\Phi^*$ to see that 
$$\bbA(a',a)\leq\Phi^*(a',-)\tensor\Phi(-,a)=\bbB(Fa',-)\tensor\bbB(-,Fa)=\bbB(Fa',Fa).$$
And for $a\in\bbA_0$ and $b\in\bbB_0$ use the regularity of $\Phi\:\bbA\dist\bbB$ (using \ref{9}) for
$$\bbB(y,Fa)\circ\bbA(a,a)=\Phi(y,a)\circ\bbA(a,a)=\Phi(y,a)=\bbB(y,Fa),$$ and likewise for $\bbA(a,a)\circ\bbB(Fa,b)=\bbB(Fa,b)$. So $F\:\bbA\to\bbB\:a\mapsto Fa$ is a regular semifunctor and, by construction, $\Phi$ converges to $F$.
\endofproof

\subsection*{Cauchy completion}

Not every totally regular $\Q$-semicategory is Cauchy complete, but any given totally regular $\Q$-semicategory $\bbB$ admits a `Cauchy completion' $\bbB\cc$ to which it is `Morita equivalent', so that finally the quantaloid $\TRSDist(\Q)$ is equivalent to its full subquantaloid determined by the Cauchy complete objects. This is what we will explain here, starting with the construction of the `Cauchy completion'.
\begin{proposition}\label{14}
For a totally regular $\Q$-semicategory $\bbB$, the following defines a totally regular $\Q$-semicategory\footnote{But $\bbB\cc$ is in general not a $\Q$-category: for a left adjoint $\phi\:*_e\dist\bbB$ in $\TRSDist(\Q)$, even though $e\leq\bbB\cc(\phi,\phi)$, it need not be the case that $1_{\dom(e)}\leq\bbB\cc(\phi,\phi)$.} $\bbB\cc$: take as objects
$$\Big\{\phi\:*_{e}\dist\bbB\ \Big|\ \mbox{$e$ is an idempotent in $\Q$, $\phi$ is a left adjoint in $\TRSDist(\Q)$}\Big\}$$
and say that the type of such an object is $t(\phi\:*_{e}\dist\bbB)=\dom(e)$; for two such objects, say $\phi\:*_e\dist\bbB$ and $\psi\:*_f\dist\bbB$, the $\Q$-arrow $\bbB\cc(\psi,\phi)\:\dom(e)\to\dom(f)$ is (the single element of) the lifting $[\psi,\phi]$ in the quantaloid $\TRSDist(\Q)$. \end{proposition}
\proof
Fist note that, by general properties of liftings, for $\phi,\psi\in\bbB\cc$, $\bbB\cc(\psi,\phi)=[\psi,\phi]=\psi^*\tensor\phi$ (where $\psi\dashv\psi^*$); this is very useful in many calculations that follow (in this proof as well as further on).
Now we verify the conditions in \ref{8} for $\bbB\cc$. For  $\phi,\psi,\theta\in\bbB\cc$ the inequality $\bbB\cc(\theta,\psi)\circ\bbB\cc(\psi,\phi)\leq\bbB\cc(\theta,\phi)$ follows trivially from the properties of liftings, whereas
$$[\psi,\phi]\tensor[\phi,\phi]
= \psi^*\tensor\phi\tensor\phi^*\tensor\phi 
= \psi^*\tensor\phi 
= [\psi,\phi]$$
implies that $\bbB\cc(\psi,\phi)\circ\bbB\cc(\phi,\phi)=\bbB\cc(\psi,\phi)$, and
likewise for $\bbB\cc(\psi,\psi)\circ\bbB\cc(\psi,\phi)=\bbB\cc(\psi,\phi)$.
\endofproof
The `stability of objects' in a totally regular $\Q$-semicategory $\bbB$ says precisely that every object $b\in\bbB_0$ is ``constant'': the left adjoint $\bbB(-,b)\:*_{\bbB(b,b)}\dist\bbB$ in $\TRSDist(\Q)$ converges to (the canonical regular semifunctor ``pointing at'') $b$---see (\ref{7.1}). Those left adjoints $\bbB(-,b)\:*_{\bbB(b,b)}\dist\bbB$ being objects of $\bbB\cc$, $\bbB$ can be now be identified within $\bbB\cc$:
$$\bbB_0\to(\bbB\cc)_0\:b\mapsto\Big(\bbB(-,b)\:*_{\bbB(b,b)}\dist\bbB\Big)$$
is a well-defined type-preserving object mapping. Denoting  $k_{\bbB}\:\bbB\to\bbB\cc$ for this map we have moreover that
$$\bbB\cc(k_{\bbB}b',k_{\bbB}b)=\bbB(b',-)\tensor\bbB(-,b)=\bbB(b',b)$$
so -- by \ref{11} -- $k_{\bbB}$ is a fully faithful regular semifunctor. It is called the {\em Cauchy completion of $\bbB\cc$}.
\par
Actually, the Cauchy completion $k_{\bbB}\:\bbB\to\bbB\cc$ of a totally regular $\Q$-semicategory is more than just a fully faithful regular semifunctor. To start with, there is a kind of ``Yoneda lemma'' for $k_{\bbB}$.
\begin{proposition}\label{15}
For any totally regular $\Q$-semicategory $\bbB$, any $b\in\bbB_0$ and any $\phi\in\bbB\cc$, $\bbB\cc(k_{\bbB}b,\phi)=\phi(b)$ and $\bbB\cc(\phi,k_{\bbB}b)=\phi^*(b)$.
\end{proposition}
\proof
Quite simply, $\bbB\cc(k_{\bbB}b,\phi)=(k_{\bbB}b)^*\tensor\phi=\bbB(b,-)\tensor\phi=\phi(b)$. The other equation is similar.
\endofproof
\begin{proposition}\label{16}
For any totally regular $\Q$-semicategory $\bbB$, the Cauchy completion $k_{\bbB}\:\bbB\to\bbB\cc$ induces an isomorphism in $\TRSDist(\Q)$: the regular semidistributors
$$\xymatrix@=15mm{
\bbB\ar@/^/@<1mm>[r]^{\bbB\cc(-,k_{\bbB}-)}|{\distsign} & \bbB\cc\ar@/^/@<1mm>[l]^{\bbB\cc(k_{\bbB}-,-)}|{\distsign}}$$
are each other's inverse.
\end{proposition}
\proof
Since $k_{\bbB}\:\bbB\to\bbB\cc$ is a fully faithful regular semifunctor, the unit inequality of the adjunction $\bbB\cc(-,k_{\bbB}-)\dashv\bbB\cc(k_{\bbB}-,-)$ is an equality: for $b,b'\in\bbB_0$,
$$\bbB\cc(k_{\bbB}b',-)\tensor\bbB\cc(-,k_{\bbB}b)=\bbB\cc(k_{\bbB}b',k_{\bbB}b)
=\bbB(b',b).$$
And for any two objects $\psi,\phi$ of $\bbB\cc$, using \ref{15},
$$\bbB\cc(\psi,k_{\bbB}-)\tensor\bbB\cc(k_{\bbB}-,\phi)
= \psi^*\tensor\phi
= \bbB\cc(\psi,\phi),$$
so the co-unit inequality is an equality too. 
\endofproof
\begin{proposition}\label{17}
For any totally regular $\Q$-semicategory $\bbB$, its Cauchy completion $\bbB\cc$ is a Cauchy complete totally regular $\Q$-semicategory.
\end{proposition}
\proof
Consider an idempotent $e$ in $\Q$ and a left adjoint $\Phi\:*_e\dist\bbB\cc$ in $\TRSDist(\Q)$. Composing $\Phi$ with the isomorphism $\bbB\cc(k_{\bbB}-,-)\:\bbB\cc\dist\bbB$ gives a new adjunction in $\TRSDist(\Q)$ like so: 
$$\xymatrix@=15mm{
\ast_e\ar@{}[r]|{\perp}\ar@/^3mm/@<1mm>[r]|{\distsign}^{\Phi} & \bbB\cc\ar@{}[r]|{\perp}\ar@/^3mm/@<1mm>[r]|{\distsign}^{\bbB\cc(k_{\bbB}-,-)}\ar@/^3mm/@<1mm>[l]|{\distsign}^{\Phi^*} & \bbB\ar@/^3mm/@<1mm>[l]|{\distsign}^{\bbB\cc(-,k_{\bbB}-)}}.$$
Denoting $\phi\dashv\phi^*$ for this adjunction, $\phi$ is an object of $\bbB\cc$ of type $t\phi=\dom(e)$ satisfying $e\leq\phi^*\tensor\phi=\bbB\cc(\phi,\phi)$. Using \ref{15} and \ref{16} we can calculate that, for any $\psi\in\bbB\cc$,
$$\bbB\cc(\psi,\phi)
= \psi^*\tensor\phi
= \bbB\cc(\psi,k_{\bbB}-)\tensor\bbB\cc(k_{\bbB}-,-)\tensor\Phi
= \bbB\cc(\psi,-)\tensor\Phi
= \Phi(\psi),$$
$$\bbB\cc(\phi,\psi)
= \phi^*\tensor\psi
=\Phi^*\tensor\bbB\cc(-,k_{\bbB}-)\tensor\bbB\cc(k_{\bbB}-,\psi)
=\Phi^*\tensor\bbB\cc(-,\psi)
=\Phi^*(\psi).$$
So $\Phi\:*_e\dist\bbB\cc$ converges to $\phi$.
\endofproof
\par
A direct consequence of all the above is the equivalence of quantaloids 
\begin{equation}\label{17.1}
\TRSDist(\Q)\simeq\TRSDist\cc(\Q)
\end{equation}
where the latter is the full subcategory of the former determined by the Cauchy complete objects.

\subsection*{Universality}

\begin{proposition}\label{18}
For two totally regular $\Q$-semicategories $\bbA$ and $\bbB$, with the latter Cauchy complete, any regular semifunctor $F\:\bbA\to\bbB$ factors in $\TRSCat(\Q)$ through the Cauchy completion $k_{\bbA}\:\bbA\to\bbA\cc$ in an essentially unique way.
\end{proposition}
\proof
Consider, for an object $\phi\:*_e\dist\bbA$ of $\bbA\cc$, the composition
$$\xymatrix@=15mm{
\ast_e\ar@{}[r]|{\perp}\ar@/^3mm/@<1mm>[r]^{\phi}|{\distsign} & 
\bbA\ar@{}[r]|{\perp}\ar@/^3mm/@<1mm>[r]^{\bbB(-,F-)}|{\distsign}\ar@/^3mm/@<1mm>[l]^{\phi^*}|{\distsign} & \bbB\ar@/^3mm/@<1mm>[l]^{\bbB(F-,-)}|{\distsign}}$$
of adjoints in $\TRSDist(\Q)$. By assumed Cauchy completeness of $\bbB$ there is, for every such $\phi$, an (essentially unique) object of $\bbB$ to which $\bbB(-,F-)\tensor\phi$ converges; call it $G\phi$. So the type of $G\phi$ is $\dom(e)=t\phi$, which already gives us a type-preserving object mapping $(\bbA\cc)_0\to\bbB_0\:\phi\mapsto G\phi$. We claim this to be the essentially unique factorization of $F$ through $k_{\bbA}$.
\par
First we verify the conditions in \ref{10} for $G$. Using that, by definition of $g\phi$, $\bbB(-,F-)\tensor\phi=\bbB(-,G\phi)$ and $\phi^*\tensor\bbB(F-,-)=\bbB(G\phi,-)$, we have for $\phi,\psi\in\bbA\cc$,
\begin{eqnarray*}
\bbA\cc(\psi,\phi)
 & = & \psi^*\tensor\phi \\
 & = & \psi^*\tensor\bbA(-,-)\tensor\phi \\
 & \leq & \psi^*\tensor\bbB(F-,F-)\tensor\phi \\
 & = & \psi^*\tensor\bbB(F-,-)\tensor\bbB(-,F-)\tensor\phi \\
 & = & \bbB(G\psi,-)\tensor\bbB(-,G\phi) \\
 & = & \bbB(G\psi,G\phi).
\end{eqnarray*}
And similarly for $\phi\in\bbA\cc$ and $b\in\bbB$,
\begin{eqnarray*}
\bbB(b,G\phi)\circ\bbA\cc(\phi,\phi)
 & = & \Big(\bbB(b,F-)\tensor\phi\Big)\circ\Big(\phi^*\tensor\phi\Big) \\
 & = & \bbB(b,F-)\tensor\Big(\phi\tensor\phi^*\tensor\phi\Big) \\
 & = & \bbB(b,F-)\tensor\phi \\
 & = & \bbB(b,G\phi),
\end{eqnarray*}
and likewise for $\bbA\cc(\phi,\phi)\circ\bbB(G\phi,b)=\bbB(G\phi,b)$. So $G$ is a regular semifunctor.
\par
Next, for any $a\in\bbA_0$, the construction of $G$ and the regularity of $F$ imply that 
$$\bbB(-,G(k_{\bbA}a))
= \bbB(-,F-)\tensor\bbA(-,a)
= \bbB(-,Fa),$$ 
proving that $\bbB(-,G\circ k_{\bbA}-)=\bbB(-,F-)$: so $G\circ k_{\bbA}\cong F$
in $\TRSCat(\Q)(\bbA,\bbB)$, that is to say, $G$ is a factorization of $F$ through $k_{\bbA}$.
\par
Finally, suppose that $H\:\bbA\cc\to\bbB$ is another factorization of $F$ through $k_{\bbA}$. Then
\begin{eqnarray*}
\bbB(-,H-)\tensor\bbA\cc(-,k_{\bbA}-)
 & = & \bbB(-,H\circ k_{\bbA}-) \\
 & = & \bbB(-,F-) \\
 & = & \bbB(-,G\circ k_{\bbA}-) \\
 & = & \bbB(-,G-)\tensor\bbA\cc(-,k_{\bbA}-)
\end{eqnarray*}
which implies that $\bbB(-,G-)=\bbB(-,H-)$ since $\bbA\cc(-,k_{\bbA}-)\:\bbA\dist\bbA\cc$ is an isomorphism. So $H\cong G$ in $\TRSCat(\Q)(\bbA\cc,\bbB)$, proving the essential uniqueness of the factorization of $F$ through $k_{\bbA}$.
\endofproof
In other words, any Cauchy completion $k_{\bbA}\:\bbA\to\bbA\cc$ of a totally regular $\Q$-semi\-cat\-e\-go\-ry $\bbA$ determines, for any Cauchy complete totally regular $\Q$-semicategory $\bbB$, a (natural) equivalence of orders
$$\TRSCat\cc(\Q)(\bbA\cc,\bbB)\iso\TRSCat(\Q)(\bbA,\bbB)\: F\mapsto F\circ k_{\bbA}.$$
So there is an adjunction between 2-categories
\begin{equation}\label{18.1}
\TRSCat\cc(\Q)\xymatrix@=15mm{
\ar@/_3mm/@<-1mm>[r]_{i}\ar@{}[r]|{\perp} & \ar@/_3mm/@<-1mm>[l]_{(-)\cc}}\TRSCat(\Q),
\end{equation}
where the former is the full sub-2-category of latter defined by the Cauchy complete objects.
\par
The point is now that the locally ordered functor
$$\TRSCat\cc(\Q)\to\Map(\TRSDist\cc(\Q))\:\Big(F\:\bbA\to\bbB\Big)\mapsto\Big(\bbB(-,F-)\:\bbA\dist\bbB\Big)$$
is the identity on objects and locally an equivalence: it is thus a biequivalence. Hence, recalling that $\TRSDist(\Q)$ is equivalent to its full subquantaloid of Cauchy complete objects, we may record a biequivalence of locally ordered categories:
$$\TRSCat\cc(\Q)\simeq\Map(\TRSDist\cc(\Q))\simeq\Map(\TRSDist(\Q)).$$
This, plus the fact that $\TRSCat\cc(\Q)$ is  full in $\TRSCat(\Q)$, in turn imply the following.
\begin{proposition}\label{19}
For totally regular $\Q$-semicategories $\bbA$ and $\bbB$, the following are equivalent:
\begin{enumerate}
\item $\bbA\cong\bbB$ in $\TRSDist(\Q)$;
\item $\bbA\cc\simeq\bbB\cc$ in $\TRSCat(\Q)$\footnote{The equivalence here must be understood in the 2-categorical sense: there exist arrows $F\:\bbA\cc\to\bbB\cc$ and $G\:\bbB\cc\to\bbA\cc$ in the locally ordered category $\TRSCat(\Q)$ satisfying $G\circ F\cong 1_{\bbA\cc}$ and $F\circ G\cong 1_{\bbB\cc}$.}.
\end{enumerate}
\end{proposition}
In words, this says that ``totally regular $\Q$-semicategories are {\em Morita equivalent} if and only if their Cauchy completions are equivalent''.

\subsection*{Comments and warnings}

{\bf (1)} The construction made in \ref{14} makes perfectly sense for any regular $\Q$-semi\-cat\-e\-go\-ry $\bbB$ (then working in the quantaloid $\RSDist(\Q)$ rather than $\TRSDist(\Q)$ of course). However, the obtained structure $\bbB\cc$ is {\em always} a totally regular $\Q$-semicategory (the proof of \ref{14} may be copied word by word). So if one wants to have a fully faithful embedding of $\bbB$ into $\bbB\cc$, then by \ref{11} necessarily $\bbB$ must be totally regular too. This explains why a ``suitable'' theory of Cauchy completions for $\Q$-semicategories only makes sense for totally regular ones (and not for merely regular ones).
\par\medskip\noindent
{\bf (2)} Because every $\Q$-category $\bbC$ is also a totally regular $\Q$-semicategory, two notions of ``Cauchy completeness'' apply: on the one hand, $\bbC$ may be Cauchy complete {\em as category} (as recalled in section \ref{B}), on the other hand, $\bbC$ may be Cauchy complete {\em as totally regular semicategory} (as explained in section \ref{D}). These notions are very different: the latter implies the former, but is strictly stronger. (Consider $\3$, the 3-element chain, viewed as quantaloid with only one object: the only one-object $\3$-category is Cauchy complete as category, but not as totally regular semicategory.) In other words, for a $\Q$-category $\bbC$, its Cauchy completion {\em as category} is very different from its Cauchy completion {\em as totally regular semicategory}. (It is maybe unfortunate that we do not notationally distinguish between both ways of Cauchy completing a $\Q$-category. The context is supposed to make clear which of both is meant.)
\par\medskip\noindent
{\bf (3)} It has been argued in [Stubbe, 2004b] that, for a regular $\Q$-semicategory $\bbA$, the pertinent notion of ``(contravariant) presheaf on $\bbA$'' is that of the so-called {\em regular (contravariant) presheaf}: such is precisely a regular semidistributor into $\bbA$, whose domain is a one-object $\Q$-category. The regular contravariant presheaves on $\bbA$ organize themselves in a $\Q$-category $\R\bbA$, and it then turns out that two regular $\Q$-semicategories $\bbA$ and $\bbB$ are isomorphic in $\RSDist(\Q)$ if and only if the $\Q$-categories $\R\bbA$ and $\R\bbB$ are equivalent in $\Cat(\Q)$. All this applies of course to totally regular $\Q$-semicategories, and together with \ref{19} we may state that: for two totally regular $\Q$-semicategories $\bbA$ and $\bbB$, $\bbA\cong\bbB$ in $\TRSDist(\Q)$ if and only if $\R\bbA\simeq\R\bbB$ in $\Cat(\Q)$, if and only if $\bbA\cc\simeq\bbB\cc$ in $\TRSCat(\Q)$. So this `Morita equivalence' of totally regular $\Q$-semicategories formally resembles the results concerning $\Q$-categories (see [Stubbe, 2004a] for example):
for two $\Q$-categories $\bbC$ and $\bbD$, $\bbC\cong\bbD$ in $\Dist(\Q)$ if and only if $\P\bbC\simeq\P\bbD$ in $\Cat(\Q)$ (equivalent categories of presheaves), if and only if $\bbC\cc\simeq\bbD\cc$ in $\Cat(\Q)$ (Cauchy completion as category). There is one big difference though: for a $\Q$-category $\bbC$, its Cauchy completion (as category) $\bbC\cc$ is a full subcategory of its presheaf category $\P\bbC$, whereas for a totally regular $\Q$-semicategory $\bbA$, its Cauchy completion $\bbA\cc$ is not at all a substructure of its regular-presheaf category $\R\bbA$.

\section{Orders and ideals over a base quantaloid}\label{F}

It is well-known that the topos $\Sh(\Omega)$ of sheaves on a locale $\Omega$ can be described in terms of $\Omega$-sets. In [Borceux and Cruciani, 1998] a similarly elementary description of the ordered objects in $\Sh(\Omega)$ is given, in terms of so-called `complete skew $\Omega$-sets'. But thinking of a locale $\Omega$ as a one-object quantaloid, such a complete skew $\Omega$-set is precisely a totally regular $\Omega$-semicategory. The main result of Borceux and Cruciani [1998] then translates as follows: {\em For a locale $\Omega$, $\TRSCat\cc(\Omega)$ is (up to equivalence) the (locally ordered) category of ordered sheaves on $\Omega$.}
\par
Motivated by this we now propose the following.
\begin{definition}\label{26}
Given a (small) quantaloid $\Q$, a {\em $\Q$-order} is a Cauchy complete totally regular $\Q$-semicategory. An {\em ideal relation} between two $\Q$-orders is a regular semidistributor between such $\Q$-orders. And an {\em order map} between two $\Q$-orders is a regular semifunctor.
\end{definition}
\begin{theorem}\label{27}
For a small quantaloid $\Q$, the $\Q$-orders and ideal relations form a (large) quantaloid $\Idl(\Q)$, and the $\Q$-orders and order maps form a (large) locally ordered category $\Ord(\Q)$, satisfying $\Ord(\Q)\simeq\Map(\Idl(\Q))$.
\end{theorem}
From (\ref{17.1}) it follows that $\Idl(\Q)$ is equivalent to $\TRSDist(\Q)$; so when working with $\Q$-orders ``up to Morita equivalence'', one can forget about Cauchy completeness. Further, (\ref{18.1}) may be understood as a ``sheafification''.
\begin{example}\label{28}
For a locale $\Omega$, $\Ord(\Omega)$ is the category of ordered objects and order maps in the topos of sheaves on $\Omega$. 
\end{example}
Applying this to the Boolean algebra $\2$, $\Ord(\2)$ is the category of ordered sets and order maps. But in fact $\Ord(\2)\simeq\Cat(\2)$ because ``accidentally'' every totally regular $\2$-semicategory is Morita equivalent to a $\2$-category. Namely, let $\bbA$ be an object of $\TRSDist(\2)$, and suppose that $a\in\bbA_0$ is such that $1\not\leq\bbA(a,a)$. Then by total regularity of $\bbA$ it follows that, for all $a'\in\bbA_0$, $\bbA(a',a)=\bbA(a,a')=0$, i.e.~such an object is ``isolated''. Consider now the full subgraph of $\bbA$ determined by the non-isolated objects; denote $i\:\wh{\bbA}\to\bbA$ for that full embedding, which is necessarily a regular semifunctor between totally regular $\2$-semicategories (by \ref{11}). But $\wh{\bbA}$ is -- by construction -- really a $\2$-category, and the regular semidistributors $\bbA(-,i-)$ and $\bbA(i-,-)$ turn out to be each other's inverse in $\TRSDist(\2)$. So the full embedding $\Dist(\2)\to\TRSDist(\2)$ is an equivalence. Recalling that $\Cat\cc(\2)=\Cat(\2)$, it follows that $\Ord(\2)\simeq\Cat(\2)$.
\par
A similar remark applies to (totally regular semi)categories enriched in the symmetric quantale $[0, \infty]$ of extended non-negative reals, with the opposite order and quantale multiplication given by addition. Categories enriched in $[0, \infty]$ are `generalized metric spaces' (and functors are distance decreasing applications), and the categorical Cauchy completion corresponds to the metric Cauchy completion [Lawvere, 1973]. Here too any totally regular $[0, \infty]$-semicategory is Morita equivalent to its subcategory of ``non-isolated'' objects, so here too $\Dist([0, \infty])\simeq\TRSDist([0, \infty])$. As a result, Cauchy complete generalized metric spaces are $[0, \infty]$-orders.
\begin{example}\label{29}
For the quantale $[0,\infty]$, $\Ord([0, \infty])$ is the category of Cauchy complete generalized metric spaces and distance decreasing maps.
\end{example}

\section{Change of base}\label{E}

Any quantaloid $\Q$ determines a new quantaloid $\Idm(\Q)$ whose objects are the idempotent arrows in $\Q$, and in which an arrow from an idempotent $e\:A\to A$ to an idempotent $f\:B\to B$ is a $\Q$-arrow $b\:A\to B$ satisfying $b\circ e=b=f\circ b$. Composition in $\Idm(\Q)$ is done as in $\Q$, the identity in $\Idm(\Q)$ on some idempotent $e\:A\to A$ is $e$ itself, and the local order in $\Idm(\Q)$ is again that of $\Q$. Actually, the quantaloid $\Idm(\Q)$ is the universal split-idempotent completion of $\Q$ in the (illegitimate) category $\QUANT$ of quantaloids. For details and references we refer to the appendix of [Stubbe, 2004b].
\par
With this in mind, the characterizations in \ref{8} and \ref{9} of `totally regular $\Q$-semicategory' and `regular semidistributor' say that:
\begin{itemize}
\item for a totally regular $\Q$-semicategory $\bbA$, every endo-hom-arrow $\bbA(a,a)\:ta\to ta$ is an object of $\Idm(\Q)$, and every $\bbA(a',a)\:ta\to ta'$ is an arrow in $\Idm(\Q)$ from $\bbA(a,a)$ to $\bbA(a',a')$,
\item for a regular semidistributor $\Phi\:\bbA\to\bbB$ from one totally regular $\Q$-semicategory to another, every element $\Phi(b,a)\:ta\to tb$ is an arrow in $\Idm(\Q)$ from $\bbA(a,a)$ to $\bbB(b,b)$.
\end{itemize}
This indicates how we can ``reshuffle'' a totally regular $\Q$-semicategory $\bbA$ to obtain a  $\Idm(\Q)$-category $\wh{\bbA}$: its objects are the same as those of $\bbA$, $(\wh{\bbA})_0=\bbA_0$, but they have different types, $\wh{t}a=\bbA(a,a)$, and its hom-arrows are $\wh{\bbA}(a',a)=\bbA(a',a)$. And similarly, the regular semidistributor $\Phi\:\bbA\dist\bbB$ determines a distributor $\wh{\Phi}:\wh{\bbA}\dist\wh{\bbB}$ between $\Idm(\Q)$-categories: simply put $\wh{\Phi}(b,a)=\Phi(b,a)$. 
\par
Note that the endo-hom-arrows of the $\Idm(\Q)$-category $\wh{\bbA}$ are identities: for an object $a\in\wh{\bbA}$, its type is $\wh{t}a=\bbA(a,a)$, and the hom-arrow on $a$ is $\wh{\bbA}(a,a)=\bbA(a,a)$, which is the identity  (in $\Idm(\Q)$!) on the idempotent $\bbA(a,a)$. For convenience, let us say that a $\Q$-category whose endo-hom-arrows are identities, is {\em normal}. And let us denote $\Dist\norm(\Q)$ for the full subquantaloid of $\Dist(\Q)$ whose objects are those normal $\Q$-categories. It is then quite straightforward to verify that 
\begin{equation}\label{21.0}
(\wh{-}):\TRSDist(\Q)\to\Dist\norm(\Idm(\Q)):
\Big(\Phi\:\bbA\dist\bbB\Big)\mapsto
\Big(\wh{\Phi}\:\wh{\bbA}\dist\wh{\bbB}\Big)
\end{equation}
is a homomorphism of quantaloids. But there's more.
\begin{lemma}\label{21}
The homomorphism in (\ref{21.0})
is an equivalence of quantaloids.
\end{lemma}
\proof
It is easy to see that $(\wh{-}):\TRSDist(\Q)\to\Dist\norm(\Idm(\Q))$ (preserves and) reflects the order on the arrows---basically because distributors are ordered elementwise, and the local order in $\Idm(\Q)$ is inherited from $\Q$. Next, given any totally regular $\Q$-semicategories $\bbA$ and $\bbB$, a distributor $\Theta\:\wh{\bbA}\dist\wh{\bbB}$ between (normal) $\Idm(\Q)$-categories consists of arrows $\Theta(b,a)\:\bbA(a,a)\to\bbB(b,b)$ in $\Idm(\Q)$ satisfying some conditions. So, really, each $\Theta(b,a)$ is a $\Q$-arrow $\Theta(b,a)\:ta\to tb$ which is compatible with the idempotents $\bbA(a,a)\:ta\to ta$ and $\bbB(b,b)\:tb\to tb$. It can now straightforwardly be verified that $\ol{\Theta}\:\bbA\dist\bbB$, defined as $\ol{\Theta}(b,a)=\Theta(b,a)$ (an equality of $\Q$-arrows), produces an arrow in $\TRSDist(\Q)$ whose image by $(\wh{-})$ is $\Theta$. So $(\wh{-}):\TRSDist(\Q)\to\Dist\norm(\Idm(\Q))$ is already fully faithful.
\par
Let $\bbC$ be a normal $\Idm(\Q)$-category: the type of an object $c\in\bbC_0$ is an idempotent $tc\:X_c\to X_c$ in $\Q$, each endo-hom-arrow $\bbC(c,c)\:tc\to tc$ is an identity in $\Idm(\Q)$, thus $\bbC(c,c)=tc$ in $\Q$, and each hom-arrow $\bbC(c',c)\:tc\to tc'$ is really a $\Q$-arrow $\bbC(c',c)\:X_c\to X_{c'}$ that is compatible with the idempotent $\Q$-arrows $tc=\bbC(c,c)$ and $tc'=\bbC(c',c')$. It can now be verified that there is a totally regular $\Q$-semicategory $\ol{\bbC}$ with the same objects as $\bbC$, $\ol{\bbC}_0=\bbC_0$, but different types, $\ol{t}c=X_c$, and hom-arrows $\ol{\bbC}(c',c)=\bbC(c',c)$ (this is an equality of $\Q$-arrows). Moreover the image by $(\wh{-})$ of $\ol{\bbC}$ is $\bbC$, so $(\wh{-}):\TRSDist(\Q)\to\Dist\norm(\Idm(\Q))$ is surjective on objects.
\endofproof
\par
Recall that a {\em monad} $t\:A\to A$ in a quantaloid $\Q$ is an endo-arrow satisfying $1_A\leq t$ and $t\circ t\leq t$; and that this monad {\em splits} when there exists an object $B$ and arrows $f\:A\to B$, $u\:B\to A$ in $\Q$ such that $u\circ f=t$ and $f\circ u=1_B$. In this case $f\dashv u$ in $\Q$, and -- in a technically rigorous sense -- $B$ is the {\em object of (free) $t$-algebras}. For more details and the relevant references we refer to the appendix of [Stubbe, 2004a]. For a $\Q$-category $\bbA$, every endo-hom-arrow $\bbA(a,a)\:ta\to ta$ is a monad\footnote{Actually, $\bbA$ itself is a monad in $\Matr(\Q)$.} in $\Q$: $1_{ta}\leq\bbA(a,a)$ is a unit-inequality, and $\bbA(a,a)\circ\bbA(a,a)\leq\bbA(a,a)$ a composition-inequality, for the $\Q$-category $\bbA$.
\begin{lemma}\label{20}
Suppose that $\bbA$ is a $\Q$-category whose endo-hom-arrows split. Then $\bbA$ is isomorphic in $\Dist(\Q)$ to a $\Q$-category whose endo-hom-arrows are identities.
\end{lemma}
\proof
Choose, for every $a\in\bbA_0$, a splitting for $\bbA(a,a)\:ta\to ta$ in $\Q$, say
$$\xymatrix@=15mm{
ta\ar@(u,l)_{\bbA(a,a)}\ar@{}[r]|{\perp}\ar@<1mm>@/^3mm/[r]^{f_a} & A_a\ar@<1mm>@/^3mm/[l]^{u_a}}.$$
Now we claim that there is a $\Q$-category $\wt{\bbA}$ with objects $(\wt{\bbA})_0=\bbA_0$ and types $\tilde{t}a=A_a$, and hom-arrows $\wt{\bbA}(a',a)=f_{a'}\circ\bbA(a',a)\circ u_a$. And that there are distributors $\Phi:\bbA\dist\wt{\bbA}$ with $\Phi(a',a)=f_{a'}\circ\bbA(a',a)$ and $\Psi:\wt{\bbA}\dist\bbA$ with $\Psi(a',a)=\bbA(a',a)\circ u_{a}$ satisfying $\Psi=\Phi\inv$ in $\Dist(\Q)$. The verifications are long but straightforward. Note in particular that the endo-hom-arrows of $\wt{\bbA}$ are identities and that $\bbA\cong\wt{\bbA}$ in $\Dist(\Q)$.
\endofproof
It follows directly from the above that, if all monads in $\Q$ split, then the full embedding $\Dist\norm(\Q)\to\Dist(\Q)$ is an equivalence of quantaloids. In particular, for any quantaloid $\Q$ the full embedding
\begin{equation}\label{20.1}
\Dist\norm(\Idm(\Q))\to\Dist(\Idm(\Q))
\end{equation}
is an equivalence: monads in a quantaloid are always idempotent and all idempotents split in $\Idm(\Q)$, hence all monads split in $\Idm(\Q)$.
\par
So from (\ref{20.1}) and \ref{21} we can now conclude.
\begin{proposition}\label{23}
For any small quantaloid $\Q$, the (large) quantaloids $\TRSDist(\Q)$ and $\Dist(\Idm(\Q))$ are equivalent.
\end{proposition}
We may read this as: ``the calculus of totally regular semicategories and semidistributors {\em enriched in $\Q$} is equivalent to the calculus of categories and distributors {\em enriched in the split-idempotent completion of $\Q$}''. 
\par
With the theory of Cauchy completions for totally regular semicategories, respectively categories, we may further say that
$\TRSDist\cc(\Q)$ and $\Dist\cc(\Idm(\Q))$ are equivalent quantaloids, and that 
$\TRSCat\cc(\Q)$ and $\Cat\cc(\Idm(\Q))$ are
biequivalent locally ordered categories. In other words we have the following.
\begin{theorem}\label{27.1}
For a small quantaloid $\Q$, the quantaloids $\Idl(\Q)$ and $\Dist\cc(\Idm(\Q)$ are equivalent, and the locally ordered categories $\Ord(\Q)$ and $\Cat\cc(\Idm(\Q))$ are biequivalent.
\end{theorem}
\par
In [Walters, 1981] it is shown that sheaves on a locale $\Omega$ are the same thing as `symmetric, skeletal Cauchy complete $\Rel(\Omega)$-categories'. The quantaloid $\Rel(\Omega)$ has as objects the elements of $\Omega$, for two objects $u,v\in\Omega$ the hom-lattice is $\Rel(\Omega)(u,v)=\{w\in\Omega\mid w\leq u\wedge v\}$, composition is finite infimum and $u$ is the identity on $u$. Thus, since every element in a locale is idempotent, $\Rel(\Omega)$ is precisely the split-idempotent completion of $\Omega$. The following is then an implicit result of Walters [1981]: {\em For a locale $\Omega$, $\Cat\cc(\Idm(\Omega))$ is (up to equivalence) the (locally ordered) category of ordered sheaves on $\Omega$.} The result in \ref{27.1} now reconciles the quite different approaches to the subject of ``ordered sheaves on a locale'', chosen by [Borceux and Cruciani, 1998] and [Walters, 1981]---see section \ref{F}.

\end{document}